\input amstex
\magnification=\magstep1 \baselineskip=13pt
\documentstyle{amsppt}
\vsize=8.7truein \CenteredTagsOnSplits \NoRunningHeads
\def\today{\ifcase\month\or
  January\or February\or March\or April\or May\or June\or
  July\or August\or September\or October\or November\or December\fi
  \space\number\day, \number\year}

\def\PP{{\bold P}}
\def\EE{{\bold E\thinspace }}
\def\conv{\operatorname{conv}}
\def\nconv{\operatorname{nconv}}

\def\Mat{\operatorname{Mat}}
\def\diag{\operatorname{diag}}
\topmatter
\title Approximating Orthogonal Matrices by Permutation Matrices \endtitle
\author Alexander Barvinok \endauthor
\address Department of Mathematics, University of Michigan, Ann Arbor,
MI 48109-1043 \endaddress \email barvinok$\@$umich.edu
\endemail
\date October 2005 \enddate
\thanks This research was partially supported by NSF Grant DMS 0400617.
The author is grateful to Microsoft (Redmond) for hospitality
during his work on this paper.
\endthanks
\abstract Motivated in part by a problem of combinatorial
optimization and in part by analogies with quantum computations,
we consider approximations of orthogonal matrices $U$ by
``non-commutative convex combinations''$A$ of permutation matrices
of the type $A=\sum A_{\sigma} \sigma$, where $\sigma$ are
permutation matrices and $A_{\sigma}$ are positive semidefinite $n
\times n$ matrices summing up to the identity matrix. We prove
that for every $n \times n$ orthogonal matrix $U$ there is a
non-commutative convex combination $A$ of permutation matrices
which approximates $U$ entry-wise within an error of $c n^{-{1
\over 2}} \ln n$ and in the Frobenius norm within an error of $c
\ln n$. The proof uses a certain procedure of randomized rounding
of an orthogonal matrix to a permutation matrix.
\endabstract
\endtopmatter
\document

\head 1. Introduction and main results \endhead

Let $O_n$ be the orthogonal group and let $S_n$ be the symmetric
group. As is well known, $S_n$ embeds in $O_n$ by means of
permutation matrices: with a permutation $\sigma$ of $\{1, \ldots,
n\}$ we associate the $n \times n$ {\it permutation matrix}
$\pi(\sigma)$,
$$\pi_{ij}(\sigma)=\cases 1 &\text{if\ }  \sigma(j)=i \\ 0
&\text{otherwise.} \endcases$$ To simplify notation, we write
$\sigma$ instead of $\pi(\sigma)$, thus identifying a permutation
with its permutation matrix and considering $S_n$ as a subgroup of
$O_n$.

In this paper, we are interested in the following general
question:
\bigskip
$\bullet$ How well are orthogonal matrices approximated by
permutation matrices?
\bigskip
\noindent
A related question is:
\bigskip
$\bullet$ Is there a reasonable way to ``round'' an orthogonal
matrix to a permutation matrix, just like real numbers are rounded
to integers?
\bigskip

To answer the second question, we suggest a simple procedure of
{\it randomized rounding}, which, given an orthogonal matrix $U$
produces not a single permutation matrix $\sigma$ but rather a
probability distribution on the symmetric group $S_n$. Using that
procedure, we show that asymptotically, as $n \longrightarrow
+\infty$, any orthogonal matrix $U$ is approximated by a certain
non-commutative convex combination, defined below, of the
permutation matrices.

\subhead (1.1) Non-commutative convex hull \endsubhead Let $v_1,
\ldots, v_m \in V$ be vectors, where $V$ is a real vector space. A
vector
$$\split v=&\sum_{i=1}^m \lambda_i v_i \quad \text{where} \\
&\sum_{i=1}^m \lambda_i =1 \quad \text{and} \quad \lambda_i \geq 0
\quad \text{for} \quad i=1, \ldots, m \endsplit \tag1.1.1$$ is
called a {\it convex combination} of $v_1, \ldots, v_m$. The set
of all convex combinations of vectors from a given set $X \subset
V$ is called the {\it convex hull} of $X$ and denoted $\conv(X)$.
We introduce the following extension of the convex hull, which we
call the non-commutative convex hull.

Let $V$ be a Hilbert space with the scalar product $\langle \cdot,
\cdot \rangle$. Recall that a self-conjugate linear operator $A$
on $V$ is called {\it positive semidefinite} provided $\langle Av,
v \rangle \geq 0$ for all $v \in V$. To denote that $A$ is
positive semidefinite, we write $A \succeq 0$. Let $I$ denote the
identity operator on $V$.

We say that $v$ is a {\it non-commutative convex combination} of
$v_1, \ldots, v_m$ if
$$\split v=&\sum_{i=1}^m A_i v_i \quad \text{where} \\
&\sum_{i=1}^m A_i =I \quad \text{and} \quad A_i \succeq 0 \quad
\text{for} \quad i=1, \ldots, m. \endsplit \tag1.1.2$$ The set of
all non-commutative convex combinations of vectors from a given
set $X \subset V$ we call the {\it non-commutative convex hull} of
$X$ and denote $\nconv(X)$.

A result of M. Naimark \cite{Na43} describes a general way to
construct operators $A_i \succeq 0$ such that $A_1 + \ldots + A_m
=I$. Namely, let $T: V \longrightarrow W$ be an embedding of
Hilbert spaces and let $T^{\ast}: W \longrightarrow V$ be the
corresponding projection. Let $\displaystyle W=\bigoplus_{i=1}^m
L_i$ be a decomposition of $W$ into a direct sum of pairwise
orthogonal subspaces and let $P_i: W \longrightarrow L_i$ be the
orthogonal projections. We let $A_i=T^{\ast}P_iT$.

A set of non-negative numbers $\lambda_1, \ldots, \lambda_m$
summing up to $1$ can be thought of as a probability distribution
on the set $\{1, \ldots, m\}$. Similarly, a set of positive
semidefinite operators $A_i$ summing up to the identity matrix can
be thought of as a measurement in a quantum system, see, for
example, \cite{Kr05}. While we can think of a convex combination
of vectors as the expected value of a vector sampled from some set
according to some probability distribution, we can think of a
non-commutative convex combination as the expected measurement of
a set of vectors.

It is clear that $\nconv(X)$ is a convex set and that
$$\conv(X) \subset \nconv(X)$$
since we get a regular convex combination (1.1.1) if we choose
$A_i$ in (1.1.2) to be the scalar operator of multiplication by
$\lambda_i$.

\subhead (1.2) Convex hulls of the symmetric group and of the
orthogonal group
\endsubhead

The convex hull of the permutation matrices $\sigma \in S_n$,
described by the Birkhoff- von Neumann Theorem, consists of the $n
\times n$ doubly stochastic matrices $A$, that is, non-negative
matrices with all row and column sums equal to 1, see, for
example, Section II.5 of \cite{Ba02}.

The convex hull of the orthogonal matrices $U \in O_n$ consists of
all the operators of norm at most 1, that is, of the operators $A:
{\Bbb R}^n \longrightarrow {\Bbb R}^n$ such that $\|A x\| \leq
\|x\|$ for all $x \in {\Bbb R}^n$, where $\| \cdot \|$ is the
Euclidean norm on ${\Bbb R}^n$, see, for example, \cite{Ha82}.

In this paper, we consider the non-commutative convex hull
$\nconv\left(S_n\right)$ of the symmetric group and show that
asymptotically, as $n \longrightarrow +\infty$, it approximates
all the orthogonal matrices. To state our main result, we consider
the following two norms on matrices: the $\ell^{\infty}$ norm
$$\|B\|_{\infty}=\max_{i,j} |\beta_{ij}|$$
and the Frobenius or $\ell^2$ norm
$$\|B\|_F =\left( \sum_{i,j=1}^n \beta^2_{ij} \right)^{1/2},$$
where $B=(\beta_{ij})$.

We prove the following result.

\proclaim{(1.3) Theorem} For every orthogonal $n \times n$ matrix
$U$ there exist positive semidefinite $n \times n$ matrices
$A_{\sigma} \succeq 0$, $\sigma \in S_n$, such that
$$\sum_{\sigma \in S_n} A_{\sigma} =I,$$
where $I$ is the $n \times n$ identity matrix, and such that for
the non-commutative convex combination
$$A =\sum_{\sigma \in S_n} A_{\sigma} \sigma$$
we have
$$\|U-A\|_{\infty} \leq c {\ln n \over \sqrt{n}}$$
and
$$\|U-A\|_F \leq c \ln n,$$
where $c$ is an absolute constant.
\endproclaim

\subhead (1.4) Discussion \endsubhead We consider $A_{\sigma}
\sigma$ as the usual product of $n \times n$ matrices. Thus the
matrix $A_{\sigma}$ acts as a linear operator
$$X \longmapsto A_{\sigma}X$$
on the space $\Mat_n$ of $n \times n$ matrices $X$. Identifying
$$\Mat_n =\underbrace{{\Bbb R}^n \oplus \cdots \oplus {\Bbb
R}^n}_{n \text{\ times}}$$ by slicing a matrix onto its columns,
we identify the action of $A_{\sigma}$ with the block-diagonal
operator
$$\left( \matrix A_{\sigma} & 0 & \ldots  &\ldots & \ldots & 0 & 0 \\ 0 & A_{\sigma}
& \ldots &\ldots &\ldots & 0 & 0  \\ \ldots & \ldots &\ldots
&\ldots &\ldots &\ldots &\ldots  \\ \ldots & \ldots &\ldots
&\ldots &\ldots &\ldots &\ldots \\
0 & 0 & \ldots & \ldots &\ldots & A_{\sigma} & 0 \\
0 & 0 & \ldots & \ldots & \ldots & 0 & A_{\sigma}
\endmatrix \right)$$
on ${\Bbb R}^n \oplus \ldots \oplus {\Bbb R}^n$.

Hence the combination $\sum_{\sigma} A_{\sigma} \sigma$ indeed
fits the definition of Section 1.1 of a non-commutative convex
combination.

Let $v=(1, \ldots, 1)$ interpreted as a column vector. Then, for
any $A=\sum_{\sigma} A_{\sigma} \sigma$ where $\sum_{\sigma}
A_{\sigma}=I$, we have $Av=v$. In particular, if $Uv \ne v$, the
matrix $U$ cannot be exactly equal to $A$, so the asymptotic
character of Theorem 1.3 is unavoidable. Taking $U=-I$ we note
that one cannot approximate $U$ entry-wise better than within
$1/n$ error, say. If $U$ is a ``typical'' orthogonal matrix, then
we have $\|U\|_{\infty} \approx c_1 \sqrt{n^{-1} \ln n}$ for some
absolute constant $c_1$, cf., for example, Chapter 5 of
\cite{MS86}. It follows from our proof that for such a typical $U$
we will have $\|U-A\|_{\infty} \leq c_2 n^{-1} \ln n$ for some
other absolute constant $c_2$.

We also note that $\|U\|_F=\sqrt{n}$ for every $U \in O_n$, so the
error in the Frobenius norm is exponentially small compared to the
norm of the matrix.

It is a legitimate question whether the bounds in Theorem 1.3 can
be sharpened.

One can ask what kind of matrices one can expect to get via
non-commutative convex combinations
$$A=\sum_{\sigma \in S_n} A_{\sigma} \sigma$$
of permutation matrices. It is easy to notice that the resulting
matrices $A$ can be quite far away from the (usual) convex hull of
the orthogonal matrices. Consider, for example, the following
situation: for the identity permutation $\sigma$, let $A_{\sigma}$
be the projection onto the first coordinate, for every
transposition $\sigma=(1k)$, $k=2, \ldots, n$, let $A_{\sigma}$ to
be the projection onto the $k$th coordinate, and for all other
$\sigma$, let $A_{\sigma}=0$. Then, for $A=(\alpha_{ij})$ we have
$\alpha_{1i}=1$ for all $i$ and all other entries of $A$ are 0.
Thus the operator norm of $A$ is $\sqrt{n}$.

\subhead (1.5) Rounding an orthogonal matrix to a permutation
matrix \endsubhead The key construction used in the proof of
Theorem 1.3 is that of a randomized rounding of an orthogonal
matrix to a permutation matrix. By now, the idea of randomized
rounding (be it the rounding of a real number to an integer or the
rounding of a positive semidefinite matrix to a vector) proved
itself to be extremely useful in optimization and other areas,
see, for example, \cite{MR95}. Let $U$ be an $n \times n$
orthogonal matrix and let $x \in {\Bbb R}^n$ be a vector. Let
$y=Ux$, so
$$x=(\xi_1, \ldots, \xi_n) \quad \text{and} \quad
y=(\eta_1, \ldots, \eta_n).$$ Suppose that the coordinates $\xi_i$
of $x$ are distinct and that the coordinates $\eta_i$ of $y$ are
distinct.
 Let
$\phi, \psi: \{1, \ldots, n \} \longrightarrow \{1, \ldots, n\}$
be the orderings of the coordinates of $x$ and $y$ respectively:
$$\xi_{\phi(1)} < \xi_{\phi(2)}< \ldots < \xi_{\phi(n)} \quad
\text{and} \quad \eta_{\psi(1)} < \eta_{\psi(2)} < \ldots <
\eta_{\psi(n)}.$$ We define the {\it rounding of} $U$ {\it at} $x$
as the permutation $\sigma=\sigma(U, x)$, $\sigma \in S_n$, such
that
$$\sigma\bigl(\phi(k)\bigr)=\psi(k) \quad \text{for} \quad
k=1, \ldots, n.$$ In words: $\sigma=\sigma(U, x)$ matches the
$k$th smallest coordinate of $x$ with the $k$th smallest
coordinate of $y=Ux$ for $k=1, \ldots, n$.

Let $\mu_n$ be the standard Gaussian measure on ${\Bbb R}^n$ with
the density
$$(2 \pi)^{-n/2} e^{-\|x\|^2/2} \quad \text{where} \quad
\|x\|^2=\xi_1^2 + \ldots + \xi_n^2 \quad \text{for} \quad
x=(\xi_1, \ldots, \xi_n).$$

If we sample $x \in {\Bbb R}^n$ at random with respect to $\mu_n$
then with probability 1 the coordinates of $x$ are distinct and
the coordinates of $y=Ux$ are distinct. Thus the rounding
$\sigma(U,x)$ is defined with probability 1. Fixing $U$ and
choosing $x$ at random, we obtain a certain probability
distribution on the symmetric group $S_n$.

The crucial observation is that for a typical $x$, the vector
$y=Ux$ is very close to the vector $\sigma x$ for
$\sigma=\sigma(U, x)$. In other words, the action of a given
orthogonal matrix on a random vector $x$ with high probability is
very close to a permutation of the coordinates. However, the
permutation varies as $x$ varies.

We prove the following result.

\proclaim{(1.6) Theorem} Let $U$ be an $n \times n$ orthogonal
matrix. For $x \in {\Bbb R}^n$, let $\sigma(U,x) \in S_n$ be the
rounding of $U$ at $x$. Let $z(x)=x-\sigma(U, x)x$ and let
$\zeta_i(x)$ be the $i$th coordinate of $z(x)$. Then
$$\int_{{\Bbb R}^n} \zeta_i^2(x) \ d \mu_n(x) \leq c {\ln^2 n \over n}$$
for some absolute constant $c$ and $i=1, \ldots, n$.
\endproclaim
\subhead (1.7) Discussion \endsubhead It follows from Theorem 1.6
that
$$\int_{{\Bbb R}^n} \|z(x)\|^2 \ d \mu_n(x) \leq c \ln^2 n.$$
Thus, for a typical $x \in {\Bbb R}^n$, we should have
$$\|Ux - \sigma(U, x) x \|=O\left( \ln n \right).$$
This should be contrasted with the fact that for a typical $x \in
{\Bbb R}^n$ we have
$$\|x\| \approx n^{1/2}.$$
Indeed, for any $0 < \epsilon <1$, we have
$$\split &\mu_n \Bigl\{x \in {\Bbb R}^n: \quad \|x\|^2 > {n \over
1-\epsilon} \Bigr\} \leq \exp\left\{-{\epsilon^2 n \over 4}
\right\} \quad
\text{and} \\
&\mu_n \Bigl\{ x \in {\Bbb R}^n: \quad \|x\|^2 \leq (1-\epsilon) n
\Bigr\} \leq \exp\left\{-{\epsilon^2 n \over 4} \right\},
\endsplit$$
see, for example, Section V.5 of \cite{Ba02}.

 Thus, for on a typical $x$, the action of
operator $U$ and the permutation $\sigma(U,x)$ do not differ much.
\bigskip
The paper is structured as follows.

In Section 2, we discuss some general properties of the proposed
randomized rounding and its possible application in the Quadratic
Assignment Problem, a hard problem of combinatorial optimization.

In Section 3, we establish concentration inequalities for the
order statistics of the Gaussian distribution on which the proof
of Theorem 1.6 is based.

In Section 4, we prove Theorem 1.6.

In Section 5, we deduce Theorem 1.3 from Theorem 1.6.

In Section 6, we conclude with some general remarks.

\head 2. Randomized rounding \endhead

The procedure described in Section 1.5 satisfies some
straightforward properties that one expects a rounding procedure
to satisfy. Given a matrix $U \in O_n$, the rounding $\sigma(U,x)
\in S_n$ for $x \in {\Bbb R}^n$ is well-defined with probability
1. Thus as $x$ ranges over ${\Bbb R}^n$, with every orthogonal
matrix $U$ we associate a probability distribution $p_U$ on $S_n$:
$$p_U(\sigma)=\mu_n \Bigl\{x \in {\Bbb R}^n: \quad \sigma(U,
x)=\sigma \Bigr\}.$$ In other words, $p_U(\sigma)$ tells us how
often do we get a particular permutation $\sigma \in S_n$ as a
rounding of $U$. For example, if $U=-I$ then $p_U$ is uniform on
the permutations $\sigma$ that are the products of $\lfloor n/2
\rfloor$ commuting transpositions: $\sigma(-I, x)$ is the
permutation matching the smallest $k$th coordinate of $x$ to its
$(n-k)$th smallest coordinate.

We note that if $U$ is a permutation matrix itself, then
$\sigma(U, x)=U$ with probability 1, so permutation matrices are
rounded to themselves. By continuity, if $U$ is close to a
permutation matrix, one can expect that the distribution $p_U$
concentrates around that permutation matrix. One can also show
that if $U$ is ``local'', that is, acts on some set $J$ of $k \ll
n$ coordinates of $x$ then $\sigma(U, x)$ is also ``local'' with
high probability, that is, acts on some $s \ll n$ coordinates
containing $J$.

If $\rho \in S_n$ is a permutation then $\sigma(\rho U, x)=\rho
\sigma (U,x)$. Therefore, if we fix $x \in {\Bbb R}^n$ with
distinct coordinates and sample $U$ at random from the Haar
probability measure on $O_n$, we get a probability distribution on
$S_n$ which is invariant under the left multiplication by $S_n$
and hence is the uniform distribution. Thus, for any fixed $x \in
{\Bbb R}^n$ with distinct coordinates, the rounding of a random
matrix $U \in O_n$ is a random permutation $\sigma \in S_n$.
Geometrically, every such an $x$ produces a partition of $O_n$
onto $n!$ isometric regions, each consisting of the matrices
rounded at $x$ to a given permutation $\sigma \in S_n$.

We also note that $\sigma(U, x)=\sigma(U, -x)$.

\subhead (2.1) Rounding in the Quadratic Assignment Problem
\endsubhead
Let us define the scalar product on the space $\Mat_n$ of real $n
\times n$ matrices by
$$\big\langle A, B \big\rangle =\sum_{i, j} a_{ij} b_{ij} \quad
\text{for} \quad A=\left(a_{ij}\right) \quad \text{and} \quad
B=\left(b_{ij}\right).$$ Given two $n \times n$ matrices $A$ and
$B$, let us consider the function $f: S_n \longrightarrow {\Bbb
R}$ defined by
$$f(\sigma)=\big\langle A,\ \sigma B \sigma^{-1} \big\rangle$$
(recall that we identify $\sigma$ with its permutation matrix).
The problem of minimizing $f$ over $S_n$, known as the {\it
Quadratic Assignment Problem}, is one of the hardest combinatorial
optimization problems, see \cite{\c Ce98}. It has long been known
that if one of the matrices is symmetric (in which case the other
can be replaced by its symmetric part, so we may assume that both
$A$ and $B$ are symmetric), then an easily computable ``eigenvalue
bound'' is available. Namely, let
$$\lambda_1 \geq \lambda_2 \geq \ldots \geq \lambda_n$$
be the eigenvalues of $A$ and let
$$\mu_1 \geq \mu_2 \geq \ldots \geq \mu_n$$
be the eigenvalues of $B$. Then the minimum value of $f$ is at
least
$$\sum_{i=1}^n \lambda_i \mu_{n-i}. \tag2.1.1$$
The bound (2.1.1) comes from extending the function $f: S_n
\longrightarrow {\Bbb R}$ to the function $f: O_n \longrightarrow
{\Bbb R}$ defined by
$$f(U)=\big\langle A,\ UBU^{\ast} \big\rangle.$$
It is then easy to compute the minimum of $f$ on $O_n$.

First, we compute $U_1$ such that $U_1 B
U_1^{\ast}=\diag\left(\mu_1, \ldots, \mu_n \right)$ is the
diagonal matrix. Next, we notice that
$$\split f(U)=&\big\langle A,\ U B U^{\ast} \big\rangle=
\big\langle A,\ (U U_1^{\ast}) U_1 B U_1^{\ast} (U_1 U^{\ast}) \big\rangle\\
= &\big\langle U_1 U^{\ast} A U U_1^{\ast},\ U_1 B U_1^{\ast}
\big\rangle.
\endsplit $$
It is then easy to see that the minimum of $f(U)$ is achieved when
$U_1 U^{\ast}=U_2$ such that $U_2A
U_2^{\ast}=\diag\left(\lambda_n, \ldots, \lambda_1 \right)$. Then
we compute $U=U_2^{\ast} U_1$.

The eigenvalue bound (2.1.1) may be far off the minimum of $f$ on
$S_n$, in which case one would expect the optimal matrix $U \in
O_n$ to be far away from a single permutation matrix. Suppose, for
example, that $n=2m$ is even. Let $J$ be the $m \times m$ matrix
of all 1's and let
$$A=\left( \matrix 1 & 1 \\ 1 & 1 \endmatrix \right) \otimes J \quad
\text{and} \quad B=\left( \matrix 1 & 0 \\ 0 & -1 \endmatrix
\right) \otimes J.$$ Then $f(\sigma) \equiv 0$ on $S_n$ while the
values of $f$ on $O_n$ range from $-n^2/2$ to $n^2/2$.

However, if $U$ is close to a particular permutation matrix, that
matrix may be recovered by rounding.

\head 3. Concentration for order statistics \endhead

Let $\xi_1, \ldots, \xi_n$ be independent identically distributed
real valued random variables. We define their {\it order
statistics} as the random variables $\omega_1, \ldots, \omega_n$,
$\omega_k=\omega_k(\xi_1, \ldots, \xi_n)$ such that
$$\omega_k(\xi_1, \ldots, \xi_n)= \quad \text{the $k$th smallest among}
\quad \xi_1, \ldots, \xi_n.$$ Thus $\omega_1$ is the smallest
among $\xi_1, \ldots, \xi_n$ and $\omega_n$ is the largest among
$\xi_1, \ldots, \xi_n$. We have
$$\omega_1 \leq \omega_2 \leq \ldots \leq \omega_n.$$

We need some concentration inequalities for order statistics.
\proclaim{(3.1) Lemma} Suppose that the cumulative distribution
function $F$ of $\xi_i$ is continuous and strictly increasing. Let
$k$ be an integer, $1 \leq k \leq n$.
\roster
\item
Let $\alpha$ be a number such that $F(\alpha) < k/n <2F(\alpha)$.
Then
$$\PP\bigl\{ \omega_k <\alpha \bigr\} \leq
\exp\left\{-{n \over 3 F(\alpha)} \left({k \over n}
-F(\alpha)\right)^2 \right\}. $$
\item
Let $\alpha$ be a number such that $F(\alpha) > k/n$. Then
$$\PP\bigl\{ \omega_k >\alpha \bigr\} \leq
\exp\left\{-{n \over 2 F(\alpha)} \left({k \over n}
-F(\alpha)\right)^2 \right\}.$$
\endroster
\endproclaim
\demo{Proof} Let us define random variables $\chi_1, \ldots,
\chi_n$ by
$$\chi_i =\cases 1 &\text{if \ } \xi_i < \alpha \\ 0 &
\text{otherwise} \endcases$$ and let $\chi=\chi_1 + \ldots +
\chi_n$.

Thus $\chi_i$ are independent random variables and
$$\PP\Bigl\{\chi_i=1 \Bigr\}=F(\alpha)=p.$$
 We note that $\omega_k <\alpha $ if and only if $\chi \geq k$.
By Chernoff's inequality (see, for example, \cite{Mc89} or
\cite{Bo91}) we get for $0 < \epsilon < 1$
$$\PP\Bigl\{\chi \geq pn(1+\epsilon) \Bigr\} \leq
\exp\left\{-{\epsilon^2 p n \over 3} \right\}.$$ Choosing
$$\epsilon={k \over pn}-1={k \over F(\alpha) n}-1$$
we complete the proof in Part (1).

Similarly in Part (2), we have $\omega_k \geq  \alpha$ if and only
if $\chi \leq k-1$. By Chernoff's inequality we get for $0<
\epsilon <1$
$$\PP\Bigl\{ \chi \leq pn(1-\epsilon) \Bigr\} \leq
\exp\left\{-{\epsilon^2 pn \over 2} \right\}.$$ Choosing
$$\epsilon=1-{k \over pn}=1-{k \over F(\alpha) n},$$
we complete the proof of Part (2).
{\hfill \hfill \hfill} \qed
\enddemo

\proclaim{(3.2) Corollary}
\roster
\item Let $k$ be an integer, $1 \leq k \leq n$.
For $0 < \epsilon < 1/2$, let us define the number
$\alpha^-=\alpha^-(k, \epsilon)$ from the equation
$$F\left( \alpha^- \right)={(1-\epsilon) k \over n}.$$
Then
$$\PP\bigl\{\omega_k < \alpha^- \bigr\} \leq
\exp\left\{-{\epsilon^2 k \over 3(1-\epsilon)} \right\}.$$
\item Let $1 \leq k \leq n/2$ be an integer.
For $0 < \epsilon < 1$, let us define the number
$\alpha^+=\alpha^+(k, \epsilon)$ from the equation
$$F\left(\alpha^+\right)={(1+\epsilon) k \over n}.$$
Then $$\PP\bigl\{\omega_k > \alpha^+ \bigr\} \leq \exp\left\{-{
\epsilon^2 k \over 2(1+\epsilon)} \right\}.$$
\endroster
\endproclaim

Next, we consider the case of the identically distributed standard
Gaussian random variables with the density
$$\phi(t)={1 \over \sqrt{2 \pi}} e^{-t^2/2}$$
and the cumulative distribution function
$$F(t)={1 \over \sqrt{2 \pi}} \int_{-\infty}^t e^{-\tau^2/2} \ d
\tau.$$

\proclaim{(3.3) Lemma} Let $\xi_1, \ldots, \xi_n$ be independent
standard Gaussian random variables. Let $1 \leq k \leq n/2$ be an
integer. Let $0 < \epsilon < 1/2$ be a number and let us define
numbers $\alpha^+=\alpha^+(k, \epsilon)$ and $\alpha^-=\alpha^-(k,
\epsilon)$ from the equations
$$F\left(\alpha^+\right)={(1+\epsilon) k \over
n} \quad \text{and} \quad F\left( \alpha^- \right)={(1-\epsilon) k
\over n}.$$ Then \roster
\item
$$\PP\Bigl\{ \omega_k < \alpha^- \Bigr\} \leq
\exp\left\{-{\epsilon^2 k \over 3} \right\} \quad \text{and} \quad
\PP\Bigl\{\omega_k > \alpha^+ \Bigr\} \leq \exp\left\{-{\epsilon^2
k \over 3} \right\};$$
\item
$$0 \leq \alpha^+ -\alpha^- \leq {\epsilon \sqrt{8 \pi} \over 1-\epsilon}.$$
\endroster
\endproclaim
\demo{Proof} Part (1) is immediate from Corollary 3.2. Clearly,
$\alpha^+ - \alpha^- \geq 0$. Applying Rolle's Theorem we get
$${2k \epsilon \over n} =F\left(\alpha^+\right) -
F\left(\alpha^- \right)= \left(\alpha^+ - \alpha^-\right)
\phi(t^{\ast}) \quad \text{for some} \quad \alpha^- < t^{\ast} <
\alpha^+.$$

Using the inequality
$$F(\alpha) < e^{-\alpha^2/2}=\sqrt{2 \pi} \phi(\alpha) \quad
 \text{for} \quad\alpha
\leq 0$$ (cf. also formula (4.2) below), we get
$$\phi(t) \geq {F(t) \over \sqrt{2 \pi}} \geq {1 \over \sqrt{2
\pi}} {(1-\epsilon)k  \over n} \quad \text{for} \quad \alpha^- < t
\leq 0.$$ By symmetry,
$$\phi(t) \geq {1 \over \sqrt{2 \pi}}
{(1-\epsilon)k  \over n} \quad \text{for} \quad 0 \leq t <
\alpha^+.$$

Summarizing,
$$\alpha^+ -\alpha^- = {2k \epsilon \over n \phi(t^{\ast})} \leq
{\epsilon \sqrt{8 \pi} \over 1-\epsilon}$$ and the proof of Part
(2) follows. {\hfill \hfill \hfill} \qed
\enddemo

\head 4. Proof of Theorem 1.6 \endhead

We need a technical (non-optimal) estimate.

\proclaim{(4.1) Lemma} Let $f: {\Bbb R}^n \longrightarrow {\Bbb
R}$ be a function such that $f(\lambda x)=\lambda f(x)$ for all $x
\in {\Bbb R}^n$ and all $\lambda \geq 0$. Let
$$B=\Bigl\{x \in {\Bbb R}^n: \quad \|x\| \leq n \Bigr\}$$
be a ball of radius $n$ and let $\mu_n$ be the standard Gaussian
measure on ${\Bbb R}^n$. Then there exists a constant $c$ such
that
$$\int_{{\Bbb R}^n} f^2 \ d \mu_n \leq c
\int_B f^2 \ d \mu_n$$ for all $n$.
\endproclaim
\demo{Proof} Let $S \subset {\Bbb R}^n$ be the unit sphere.
Passing to the polar coordinates, we get
$$\int_{{\Bbb R}^n} f^2 \ d \mu_n =
(2 \pi)^{-n/2} \left(\int_{S} f^2 \ d x \right) \int_0^{+\infty}
t^{n+1} e^{-t^2/2} \ dt$$ and, similarly,
$$\int_{B} f^2 \ d \mu_n =
(2 \pi)^{-n/2} \left(\int_{S} f^2 \ d x \right) \int_0^n t^{n+1}
e^{-t^2/2} \ dt.$$ Furthermore, we have
$$\int_0^{+\infty} t^{n+1} e^{-t^2/2} \ dt =2^{n/2} \Gamma\left({n+2 \over 2} \right).$$
For all sufficiently large $n$, we have
$$t^{n+1} e^{-t^2/2} \leq e^{-t^2/4} \quad \text{for all} \quad t>
n,$$ so we have
$$\int_n^{+\infty} t^{n+1} e^{-t^2/2} \ d t \leq c$$
for some constant $c$ and all $n$. The proof now follows. {\hfill
\hfill \hfill} \qed
\enddemo

Apart from Lemma 4.1, we need the estimate:
$$\mu_n\Bigl\{x=(\xi_1, \ldots, \xi_n): \quad |\xi_i| > t \Bigr\}
\leq 2e^{-t^2/2} \tag4.2$$ for any $t \geq 0$ and any $i=1,
\ldots, n$, see, for example, Section V.5 of \cite{Ba02}. Now we
can prove Theorem 1.6.

\demo{Proof of Theorem 1.6} Let $B$ be the ball of radius $n$ in
${\Bbb R}^n$ centered at the origin. By Lemma 4.1 it suffices to
prove the estimate for the integral
$$\int_B \zeta^2_i(x) \ d \mu_n(x).$$
Without loss of generality, we may assume that $i=1$ so that
$\zeta(x)=\zeta_1(x)$ is the first coordinate of
$Ux-\sigma(U,x)x$.

Let $V_k \subset B$ be the subset of $x \in B$ such the first
coordinate of $Ux$ is the $k$th smallest among the coordinates of
$Ux$. Then $V_1, \ldots, V_n$ are polyhedral (generally,
non-convex) sets that cover $B$ and intersect only at boundary
points. Since $B$ is $O_n$-invariant, the sets $V_k$ are isometric
and so we have
$$\mu_n(V_1)=\ldots = \mu_n(V_n) = {\mu_n(B) \over n} < {1 \over n}.$$

Thus we have
$$\int_{B} \zeta^2(x) \ d \mu_n(x) =\sum_{k=1}^n \int_{V_k} \zeta^2(x) \ d
\mu_n(x).$$

In what follows, $c_i$ for $i=1,2, \ldots$ denote various absolute
constants.

We note that for any $x \in B$ we have $|\zeta(x)| \leq 2n$.
Moreover, by (4.2)
$$\mu_n\Bigl\{x \in V_k: \quad |\zeta(x)| \geq c_1 \sqrt{\ln n} \Bigr\}
\leq n^{-3} \quad \text{for all sufficiently large} \quad n.$$
Therefore,
$$\int_{V_k} \zeta^2(x) \ d \mu_n(x) \leq c_2 \left(\ln n \right)
n^{-1}
 \quad \text{for all} \quad k. \tag4.3$$

For
$$36\ln n \leq k \leq n/2$$ and all sufficiently large $n$ we get a better estimate
via Lemma 3.3. Namely, let us choose
$\epsilon=\epsilon_k=3k^{-1/2} \sqrt{\ln n}$ in Lemma 3.3 and let
$\alpha_k^+$ and $\alpha_k^-$ be the corresponding bounds. It
follows that for $36 \ln n \leq k \leq n/2$ and all sufficiently
large $n$ we have
$$\mu_n \Bigl\{x \in {\Bbb R}^n: \quad \omega_k(x) \notin
[\alpha_k^-, \alpha_k^+] \Bigr\} \leq n^{-3}$$ and, similarly,
$$\mu_n \Bigl\{x \in {\Bbb R}^n: \quad \omega_k(Ux) \notin
[\alpha_k^-, \alpha_k^+] \Bigr\} \leq n^{-3},$$ where
$$0 \leq \alpha_k^+-\alpha_k^- \leq c_3
k^{-{1 \over 2}} \sqrt{\ln n}.$$

Hence
$$\mu_n \Bigl\{x \in {\Bbb R}^n: \quad |\omega_k(Ux)-\omega_k(x)| >
c_3 k^{-{1 \over 2}} \sqrt{\ln n} \Bigr\} \leq 2 n^{-3}.$$ Since
for $x \in V_k$ we have $\zeta(x)=\omega_k(Ux)-\omega_k(x)$, we
conclude
$$\mu_n\Bigl\{x \in V_k: \quad |\zeta(x)| > c_3 k^{-{1 \over 2}} \sqrt{\ln n} \Bigr\}
 \leq 2 n^{-3}$$
and
$$\int_{V_k} \zeta^2(x) \ d \mu_n(x) \leq c_4 n^{-1} k^{-1} \ln n
 \quad \text{for} \quad 36 \ln n  \leq k \leq n/2 \tag4.4$$
 and all sufficiently large $n$.

Summarizing (4.3) and (4.4), we get
$$\split \sum_{1 \leq k \leq n/2} &\int_{V_k} \zeta^2(x) \ d \mu_n(x) \\=
&\sum_{1 \leq k < 36 \ln n} \int_{V_k} \zeta^2(x) \ d \mu_n(x)+
\sum_{36 \ln n \leq k \leq n/2} \int_{V_k} \zeta^2(x) \ d \mu_n(x)
\\ \leq &c_5 \left( \ln n \right)^2 n^{-1}. \endsplit$$

Since by the symmetry $x \leftrightarrow -x$ we have
$$\int_{V_k} \zeta^2(x) \ d \mu_n(x) = \int_{V_{n-k}} \zeta^2(x) \ d
\mu_n(x),$$ the proof follows. {\hfill \hfill \hfill} \qed
\enddemo

\head 5. Proof of Theorem 1.3 \endhead

First, we introduce some notation.

For vectors $x=(\xi_1, \ldots, \xi_n)$ and $y=(\eta_1, \ldots,
\eta_n)$ let $x \otimes y$ be the $n \times n$ matrix with the
$(i,j)$th entry equal to $\xi_i \eta_j$.

We observe that for any $n \times n$ matrix $A$ we have
$$A(x \otimes y)=(Ax) \otimes y,$$
where the product in the left hand side we interpret as the
product of matrices and the product $Ax$ in the right hand side we
interpret as a product of a matrix and a column vector.

Let
$$\langle x, y \rangle = \sum_{i=1}^n \xi_i \eta_i \quad
\text{for} \quad x=(\xi_1, \ldots, \xi_n) \quad \text{and} \quad
y=(\eta_1, \ldots, \eta_n)$$ be the standard scalar product in
${\Bbb R}^n$. Then for all $x, y, a \in {\Bbb R}^n$, we have
$$\left(x \otimes y \right) a= \langle a, y \rangle x.\tag5.1$$

Let
$$\|x\|=\sqrt{\langle x, x \rangle} \quad \text{for} \quad x \in {\Bbb R}^n$$
be the usual Euclidean norm of a vector.

We need a couple of technical results.

\proclaim{(5.2) Lemma} Let $L$ be an $n \times n$ matrix. Then
$$\|L\|_F^2=\int_{{\Bbb R}^n} \|L a\|^2 \ d \mu_n(a).$$
\endproclaim
\demo{Proof} Let $a=(\alpha_1, \ldots, \alpha_n)$, where
$\alpha_i$ are independent standard Gaussian random variables.
Then
$$\|La\|^2=\sum_{i=1}^n \left(\sum_{j=1}^n l_{ij} \alpha_j
\right)^2.$$ Since $\EE \alpha_i \alpha_j = 0$ for $i \ne j$ and
$\EE \alpha_j^2=1$, taking the expectation we get
$$\EE \|La\|^2 = \sum_{i,j=1}^n l_{ij}^2.$$
{\hfill \hfill \hfill} \qed
\enddemo

\proclaim{(5.3) Lemma} Let $f: {\Bbb R}^n \longrightarrow {\Bbb
R}$ be an integrable function such that $$\int_{{\Bbb R}^n} f^2(x)
\ d \mu_n(x)< +\infty \quad \text{and} \quad  \int_{{\Bbb R}^n}
\|x\|^2 f^2(x) \ d \mu_n(x) < + \infty.$$ Then
$$\int_{{\Bbb R}^n} \left( \int_{{\Bbb R}^n} \langle a, x \rangle
f(x) \ d \mu_n(x) \right)^2 \ d \mu_n(a) \leq \int_{{\Bbb R}^n}
f^2(x) \ d \mu_n(x).$$
\endproclaim
\demo{Proof} Let ${\Cal L}$ be the subspace of the Hilbert space
$L^2\left({\Bbb R}^n, \mu_n \right)$ consisting of the linear
functions and let ${\Cal L}^{\bot}$ be its orthogonal complement.
We write
$$f= \langle b, x \rangle + h,$$
where $b \in {\Bbb R}^n$ and $h \in {\Cal L}^{\bot}$.

Hence we have
$$\int_{{\Bbb R}^n} f^2(x) \ d \mu_n(x) \geq \int_{{\Bbb R}^n}
\langle b, x \rangle^2 \ d \mu_n(x)= \langle b, b \rangle$$ and
$$\int_{{\Bbb R}^n} \langle a, x \rangle f(x) \ d \mu_n(x) =
\int_{{\Bbb R}^n} \langle a, x \rangle \langle b, x \rangle \ d
\mu_n(x) = \langle a, b \rangle.$$

Therefore,
$$\split &\int_{{\Bbb R}^n} \left( \int_{{\Bbb R}^n} \langle a, x \rangle
f(x) \ d \mu_n(x) \right)^2 \ d \mu_n(a)=\int_{{\Bbb R}^n} \langle
a, b \rangle^2 \ d \mu_n(a)\\ =&\langle b, b \rangle \leq
\int_{{\Bbb R}^n} f^2(x) \ d \mu_n(x) \endsplit$$ as claimed.
{\hfill \hfill \hfill} \qed
\enddemo

Now we are ready prove Theorem 1.3.

\demo{Proof of Theorem 1.3} Given an orthogonal matrix $U$, we
will construct a matrix $A$ approximating $U$ as desired in the
form
$$A=\sum_{\sigma \in S_n} \sigma A_{\sigma}, \quad \text{where} \quad
\sum_{\sigma \in S_n} A_{\sigma} =I \quad \text{and} \quad
A_{\sigma} \succeq 0.$$ To get the approximation of the type
$$A=\sum_{\sigma \in S_n} A_{\sigma} \sigma$$
claimed in the Theorem, one should apply the construction to
$U^{\ast}$.

Let $\sigma(U, x)$ be the rounding of $U$ at $x$ and let us define
$$X_{\sigma}=\Bigl\{x \in {\Bbb R}^n: \quad \sigma(U,x)=\sigma
\Bigr\} \quad \text{for} \quad \sigma \in S_n$$ and
$$A_{\sigma}=\int_{X_{\sigma}} x \otimes x \ d \mu_n(x).$$
Clearly, $A_{\sigma}$ are positive semidefinite and
$$\sum_{\sigma \in S_n} A_{\sigma}=\int_{{\Bbb R}^n} x \otimes x \ d \mu_n(x)=I.$$
On the other hand,
$$\split U - \sum_{\sigma} \sigma A_{\sigma}=
&\int_{{\Bbb R}^n} \left(U x \right) \otimes x \ d \mu_n(x)-
\sum_{\sigma \in S_n} \sigma \int_{X_{\sigma}} x \otimes x \ d
\mu_n(x) \\= &\int_{{\Bbb R}^n} \left(U x \right) \otimes x \ d
\mu_n(x)- \sum_{\sigma \in S_n} \int_{X_{\sigma}} \sigma(x)
\otimes x \ d \mu_n(x)\\=&\int_{{\Bbb R}^n} \left(Ux - \sigma(U,x)
x \right) \otimes x \ d \mu_n(x). \endsplit$$

Let
$$L=\int_{{\Bbb R}^n} \left(Ux - \sigma(U,x) x \right)
\otimes x \ d \mu_n(x)=\int_{{\Bbb R}^n} z(x) \otimes x \ d
\mu_n(x)$$ in the notation of Theorem 1.3. Thus $L$ is an $n
\times n$ matrix, $L=(l_{ij})$ and
$$U-A=L.$$
 Using Theorem 1.6, we estimate $l_{ij}$. Denoting $\xi_j(x)$ the $j$th coordinate
 of $x$, we get from Theorem 1.6
$$\split |l_{ij}|=\Big| \int_{{\Bbb R}^n} \zeta_i(x) \xi_j(x) \ &d \mu_n(x) \Big| \leq
\left(\int_{{\Bbb R}^n} \zeta^2_i(x) \ d \mu_n(x) \right)^{1/2}
\left(\int_{{\Bbb R}^n} \xi_j^2(x) \ d \mu_n(x) \right)^{1/2} \\
\leq & c n^{-{1 \over 2}} \ln n, \endsplit$$ from which we get
$$\|U-A\|_{\infty} \leq c n^{-{1 \over 2}} \ln n$$
as desired.

Finally, we estimate
$$\|U-A\|_F=\|L\|_F$$
using Lemma 5.2. By formula (5.1) for $a \in {\Bbb R}^n$ we have
$$La=\int_{{\Bbb R}^n} \langle a, x \rangle \left(Ux -\sigma(U,x)x \right)
\ d \mu_n(x)=\int_{{\Bbb R}^n} \langle a,x \rangle z(x) \ d
\mu_n(x)$$ in the notation of Theorem 1.6.
 Let us estimate
 $$\int_{{\Bbb R}^n} \|La\|^2 \ d \mu_n(a).$$ The $i$th coordinate
$\lambda_i(a)$ of $La$ is
$$\lambda_i(a)=\int_{{\Bbb R}^n} \langle a,x \rangle \zeta_i(x) \ d \mu_n(x).$$
By Lemma 5.3,
$$\split \int_{{\Bbb R}^n} \lambda^2_i(a) \ d \mu_n(a)= &\int_{{\Bbb R}^n}
\left(\int_{{\Bbb R}^n} \langle a, x \rangle \zeta_i(x) \ d
\mu_n(x)\right)^2 \ d \mu_n(a) \\ \leq &\int_{{\Bbb R}^n}
\zeta^2_i(x) \ d \mu_n(x) \leq c {\ln^2 n \over n} \endsplit$$ by
Theorem 1.3. Therefore,
$$\|L\|_F^2 =\int_{{\Bbb R}^n} \|La\|^2 \ d \mu_n(a) =
\sum_{i=1}^n \int_{{\Bbb R}^n} \lambda_i^2(a) \ d \mu_n(a) \leq c
\ln^2 n$$ as desired. {\hfill \hfill \hfill} \qed
\enddemo

\head 6. Concluding remarks \endhead

A somewhat stronger estimate follows from our proof of Theorem
1.3. Namely, let $u_1, \ldots, u_n$ be the column vectors of $U$
and let $a_1, \ldots, a_n$ be the column vectors of $A$. Then
$$\|u_i -a_i\| \leq c {\ln n \over \sqrt{n}} \quad \text{for} \quad
i=1, \ldots, n,$$ where $\| \cdot\|$ is the Euclidean norm in
${\Bbb R}^n$.

It follows from our construction of matrices $A_{\sigma}$ in the
proof of Theorem 1.3 that the trace of $A_{\sigma}$ is equal to
the probability that the matrix $U^{\ast}$ is rounded to the
permutation $\sigma$.

One can easily construct small approximate non-commutative convex
combinations
$$U \approx \sum_{i=1}^N A_i \sigma_i$$
with
$$A_i \succeq 0 \quad \text{and} \quad \sum_{i=1}^N A_i
\approx I$$ by sampling $N$ points $x_i$ at random from the
Gaussian distribution $\mu_n$, computing the rounding
$\sigma_i=\sigma\left(U^{\ast}, x_i\right)$ and letting $A_i = x_i
\otimes x_i$.


\Refs

\widestnumber\key{AAAA}

\ref\key{Ba02} \by A. Barvinok \book A Course in Convexity
\bookinfo Graduate Studies in Mathematics \vol 54 \publ American
Mathematical Society \publaddr Providence, RI \yr 2002 \endref

\ref\key{Bo91} \by B. Bollob\'as \paper Random graphs \inbook
Probabilistic combinatorics and its applications (San Francisco,
CA, 1991) \pages 1--20 \bookinfo Proc. Sympos. Appl. Math. \vol 44
\publ Amer. Math. Soc. \publaddr Providence, RI \yr 1991
\endref

\ref \key{\c Ce98} \by E. \c Cela \book The Quadratic Assignment
problem. Theory and Algorithms \bookinfo Combinatorial
Optimization \vol 1 \publ Kluwer Academic Publishers \publaddr
Dordrecht \yr 1998 \endref

\ref\key{Ha82} \by P.R. Halmos \book A Hilbert Space Problem Book.
Second edition \bookinfo Graduate Texts in Mathematics, vol. 19.
Encyclopedia of Mathematics and its Applications, vol. 17 \publ
Springer-Verlag \publaddr New York-Berlin \yr 1982 \endref

\ref \key{Kr05} \by D.W. Kribs \paper A quantum computing primer
for operator theorists \jour Linear Algebra Appl. \vol 400 \yr
2005 \pages 147--167 \endref

\ref\key{Mc89} \by C. McDiarmid \paper On the method of bounded
differences \inbook Surveys in combinatorics 1989 (Norwich, 1989)
\pages 148--188 \bookinfo London Math. Soc. Lecture Note Ser. \vol
141 \publ Cambridge Univ. Press \publaddr Cambridge \yr 1989
\endref

\ref\key{MR95} \by R. Motwani and P. Raghavan \book Randomized
Algorithms \publ Cambridge University Press \publaddr Cambridge
\yr 1995 \endref

\ref\key{MS86} \by V.D. Milman and G. Schechtman \book Asymptotic
Theory of Finite-Dimensional Normed Spaces. With an appendix by M.
Gromov \bookinfo Lecture Notes in Mathematics \vol 1200 \publ
Springer-Verlag \publaddr Berlin \yr 1986
\endref

\ref \key{Na43} \by M.A. Neumark \paper On a representation of
additive operator set functions \jour C. R. (Doklady) Acad. Sci.
URSS (N.S.) \vol 41 \yr 1943 \pages 359--361 \endref

\endRefs
\enddocument
\end